\documentclass{amsart}

\usepackage[backref=section]{hyperref}

\newtheorem{theorem}{Theorem}
\newtheorem{proposition}{Proposition}

\newcommand{\csillag}{\sideset{}{^*}}
\newcommand{\smod}[1]{\text{\rm{ (mod $#1$)}}}
\newcommand{\al}{\alpha}
\newcommand{\dd}{\delta}
\newcommand{\ep}{\epsilon}
\newcommand{\FI}{\mathfrak{I}}
\newcommand{\gG}{\gamma}
\newcommand{\GG}{\Gamma}
\newcommand{\om}{\omega}
\newcommand{\si}{\sigma}
\newcommand{\vfi}{\varphi}
\newcommand{\lL}{\lambda}

\newcommand{\lfi}{\lambda_\phi}
\newcommand{\lpsi}{\lambda_\psi}
\newcommand{\blfi}{{\bar\lambda}_\phi}
\newcommand{\HH}{\mathcal{H}}
\newcommand{\QQ}{\mathbb{Q}}
\newcommand{\CQ}{\mathcal{Q}}
\newcommand{\RR}{\mathbb{R}}
\newcommand{\CC}{\mathbb{C}}
\newcommand{\ZZ}{\mathbb{Z}}

\begin{document}

\author{Gergely Harcos}
\title{An additive problem in the Fourier coefficients of cusp forms}
\address{Department of Mathematics, Princeton University, Fine Hall,
Washington Road, Princeton, NJ 08544, USA}
\email{gharcos@math.princeton.edu}


\begin{abstract}We establish an estimate on sums of shifted
products of Fourier coefficients coming from holomorphic or Maass
cusp forms of arbitrary level and nebentypus. These sums are
analogous to the binary additive divisor sum which has been
studied extensively. As an application we derive, extending work
of Duke, Friedlander and Iwaniec, a subconvex estimate on the
critical line for $L$-functions associated to character twists of
these cusp forms.
\end{abstract}

\maketitle

\section{Introduction and statement of results}\label{sect6}

In the analytic theory of automorphic $L$-functions one often
encounters sums of the form
\begin{equation}\label{eq24}
D_f(a,b;h)=\sum_{am\pm bn=h} \lfi(m)\lpsi(n)f(am,bn),
\end{equation}
where $a$, $b$, $h$ are positive integers, $\lfi(m)$ (resp.
$\lpsi(n)$) are the normalized Fourier coefficients of a
holomorphic or Maass cusp form $\phi$ (resp. $\psi$) coming from
an automorphic representation of $GL(2)$ over $\QQ$ and $f$ is
some nice weight function on $(0,\infty)\times (0,\infty)$. These
sums have been studied extensively beginning with Selberg
\cite{Se} (see also Good \cite{G}) and are analogous to the
generalized binary additive divisor sum where $\lfi(m)$ and
$\lpsi(n)$ are replaced by values of the divisor function:
\[D_f^\tau(a,b;h)=\sum_{am\pm bn=h} \tau(m)\tau(n)f(am,bn).\]
The analogy is deeper than formal, because $\tau(n)$ appears as
the $n$-th Fourier coefficient of the modular form
$\frac{\partial}{\partial s}E(z,s)\big|_{s=1/2}$ where $E(z,s)$ is
the Eisenstein series for $SL_2(\ZZ)$. In general one tries to
deduce good estimates for these sums assuming the parameters $a$,
$b$, $h$ are of considerable size.

The binary additive divisor problem has an extensive history and
we refer the reader to \cite{DFI} for a short introduction. Let us
just mention that in the special case $a=b=1$ one can derive very
sharp results by employing the spectral theory of automorphic
forms for the full modular group \cite{Mo}. This approach is hard
to generalize for larger values of $a$, $b$ as one faces
difficulties with small Laplacian eigenvalues for the congruence
subgroup $\GG_0(ab)$ as well as uniformity issues.
The idea of Duke, Friedlander and Iwaniec
\cite{DFI} is to combine the more elementary $\dd$-method (a
variant of Kloosterman's refinement of the classical circle method)
with a Voronoi-type
summation formula for the divisor function and then apply Weil's
estimate for the individual Kloosterman sums
\[S(m,n;q)=\csillag\sum_{d\smod{q}}e_q\bigl(dm+{\bar
d}n\bigr)\]
that arise.

Assuming $a$, $b$ are coprime and the partial derivatives of the
weight function $f$ satisfy the estimate
\begin{equation}\label{eq1}
x^iy^jf^{(i,j)}(x,y)\ll_{i,j}
\left(1+\frac{x}{X}\right)^{-1}\left(1+\frac{y}{Y}\right)^{-1}P^{i+j}
\end{equation}
with some $P,X,Y\geq 1$ for all $i,j\geq 0$, they were able to deduce
\[D_f^\tau(a,b;h)=\int_0^\infty g(x,\mp x\pm h)\,dx+
O\bigl(P^{5/4}(X+Y)^{1/4}(XY)^{1/4+\ep}\bigr),\]
where the implied constant depends on $\ep$ only,
\[g(x,y)=f(x,y)\sum_{q=1}^\infty \frac{(ab,q)}{abq^2}c_q(h)
(\log x-\lL_{aq})(\log y-\lL_{bq}),\]
$c_q(h)=S(h,0;q)$ denotes Ramanujan's sum and
$\lL_{aq}$, $\lL_{bq}$ are constants given by
\[\lL_{aq}=2\gG+\log\frac{aq^2}{(a,q)^2}.\]
As was pointed out in \cite{DFI} the error term is smaller than
the main term whenever
\[ab\ll P^{-5/4}(X+Y)^{-5/4}(XY)^{3/4-\ep}.\]

The case $N=a=b=1$ of the sum $D_f(a,b;h)$ has been discussed in
detail via spectral theory by Jutila \cite{Ju1,Ju2}. This approach
is hard to generalize to $Nab>1$ because
one faces with the difficulty of uniformity and small Laplacian
eigenvalues similarly as in the additive divisor problem. A
different spectral approach was developed by Sarnak for all levels.
Using his estimates for triple products of eigenfunctions
\cite{Sa4} (see also \cite{P,BR}) he recently established quite
strong uniform bounds for $D_f(a,b;h)$ at least when the forms $\phi$
and $\psi$ are holomorphic. This method has the big advantage
of generalizing naturally to number fields \cite{Sa1,Sa3}. Our
aim here is to emphasize the Maass case and establish,
uniformly for all $h>0$, a nontrivial estimate on
$D_f(a,b;h)$ in the spirit of Duke, Friedlander and Iwaniec.
\begin{theorem}\label{Th1} Let $\lfi(m)$ (resp. $\lpsi(n)$)
be the normalized Fourier coefficients of a holomorphic or Maass
cusp form $\phi$ (resp. $\psi$) of arbitrary level and nebentypus
and suppose that $f$ satisfies (\ref{eq1}). Then for coprime $a$
and $b$ we have
\[D_f(a,b;h)\ll P^{11/10}(ab)^{-1/10}(X+Y)^{1/10}(XY)^{2/5+\ep},\]
where the implied constant depends only on $\ep$ and the forms
$\phi$, $\psi$.
\end{theorem}

See the next section for a precise definition of the notions in
the theorem. We note that the theorem supercedes the trivial upper
bound $(XY/ab)^{1/2}$ (following from Cauchy's inequality, see
Section~\ref{sect2}) whenever
\begin{equation}\label{eq17}
ab\ll P^{-11/4}(X+Y)^{-1/4}(XY)^{1/4-\ep}.
\end{equation}

As an application we prove a subconvex estimate on the critical
line for $L$-functions associated to character twists of a fixed
holomorphic or Maass cusp form $\phi$ of arbitrary level and
nebentypus. We shall assume that $\phi$ is a \emph{primitive form}, i.e.
a newform in the sense of \cite{AL,Li,ALi} normalized so that
$\lfi(1)=1$. Then $\lfi(m)$ $(m\geq 1)$ defines a
character of the corresponding Hecke algebra while
$\lfi(-m)=\pm\lfi(m)$ (with a constant sign) when $\phi$ is a
Maass form. In other words, $\phi$ defines a cuspidal automorphic
representation of $GL(2)$ over $\QQ$. The contragradient
representation corresponds to the primitive cusp form
$\tilde\phi(z)=\bar\phi(-\bar z)$ with Fourier coefficients
$\lambda_{\tilde\phi}(m)=\blfi(m)$. If $q$ is an integer prime to
the level and $\chi$ is a primitive Dirichlet character modulo $q$
then to the twisted primitive cusp form $\phi\otimes\chi$ is
attached the $L$-function
\[L(s,\phi\otimes\chi)=\sum_{m=1}^\infty\frac{\lfi(m)\chi(m)}{m^s}\]
which is absolutely convergent for $\Re s>1$ and has an Euler
product over the prime numbers. It has an analytic continuation to
an entire function, as shown by Hecke, and satisfies a functional
equation of the standard type. It follows from the
Phragm\'en--Lindel\"of convexity principle that for a fixed point
on the critical line $\Re s=1/2$ we have a bound
\[L(s,\phi\otimes\chi)\ll_{\ep,s,\phi}q^{1/2+\ep}.\]
By a subconvexity estimate we mean one which replaces the
convexity exponent $1/2$ by any smaller absolute constant. Upon
the Generalized Riemann Hypothesis we would have the Generalized
Lindel\"of Hypothesis which asserts that any positive exponent is
permissible. For the philosophy of breaking convexity in the
analytic theory of $L$-functions and its importance for arithmetic
we refer the reader to the excellent discussion by Iwaniec and
Sarnak \cite{IS}.
\begin{theorem}\label{Th2} Suppose that $\phi$ is a
primitive holomorphic or Maass cusp form of arbitrary level and
nebentypus. Let $\Re s=1/2$ and $q$ be an integer prime to the
level. If $\chi$ is a primitive Dirichlet character modulo $q$
then
\begin{equation}\label{eq19}
L(s,\phi\otimes\chi)\ll q^{1/2-1/54+\ep},
\end{equation}
where the implied constant depends only on $\ep$, $s$ and the form
$\phi$.
\end{theorem}

This estimate with exponent $1/2-1/22$ has been proved for holomorphic
forms of full level in \cite{DFI2} and the improved exponent
$1/2-7/130$ follows for holomorphic forms of arbitrary level as
a special case of Theorem 1 in \cite{Sa3}.
Duke, Friedlander and Iwaniec anticipated their
method to be extendible to more general $L$-functions of rank two,
and the present paper is indeed an extension of their work.

Combining the estimate (\ref{eq19}) at the central point $s=1/2$
with Waldspurger's theorem \cite{Wa} (see also \cite{K}) we get
the bound
\[c(q)\ll_\ep q^{1/4-1/108+\ep},\quad\text{$q$ square-free}\]
for the normalized Fourier coefficients of half-integral weight
forms of arbitrary level. Such a nontrivial bound is the key step
in the solution of the general ternary Linnik problem given by
Duke and Schulze-Pillot \cite{D,D-SP}.

The proof Theorem~\ref{Th1} is presented in Sections~\ref{sect1}
through \ref{sect2} and closely follows \cite{DFI}. The heart of
the argument is again a Voronoi-type formula (see
Section~\ref{sect1}) for transforming certain exponential sums
defined by the coefficients $\lfi(m)$ and $\lpsi(n)$ but this time
the level of the forms imposes some restriction on the frequencies
in the formula. As the $\dd$-method uses information at all
frequencies (and in this sense it corresponds to the classical
Farey dissection of the unit circle) we replace it (in
Section~\ref{sect3}) with another variant of the circle method
(given by Jutila \cite{Ju3}) which is more flexible in the choice
of frequencies. After the transformations we shall encounter
twisted Kloosterman sums
\[S_\chi(m,n;q)=\csillag\sum_{d\smod{q}}
\chi(d)e_q\bigl(dm+{\bar d}n\bigr),\] where $\chi$ is a Dirichlet
character mod $q$. We shall make use of the usual Weil--Estermann
bound \begin{equation}\label{eq10}\bigl|S_\chi(m,n;q)\bigr|\leq
(m,n,q)^{1/2}q^{1/2}\tau(q)\end{equation} which holds true for
these sums as well (the original proofs \cite{W,E} carry over with
minor modifications).

Sections~\ref{sect4} through \ref{sect5} are devoted to the proof
of Theorem~\ref{Th2}. In Section~\ref{sect4} we reduce
(\ref{eq19}), via an approximate functional equation, to an
inequality about certain finite sums involving about $q$ terms. In
order to prove this inequality we use the amplification method
which was introduced in \cite{DFI2}. The idea is to consider a
suitably weighted second moment of the finite sums arising from the family
$\phi\otimes\chi$ of cusp forms ($\chi$ varies, $\phi$ is fixed).
We choose the weights (called amplifiers) in such a way that one of the characters
$\chi$ is emphasized while the second moment average is still of moderate size.
This forces, by positivity, $L(s,\phi\otimes\chi)$ to be small. In
the course of evaluating the amplified second moment we
encounter diagonal and off-diagonal terms and it is the
off-diagonal contribution where Theorem~\ref{Th1} enters.

\medskip \noindent {\bf Acknowledgements.} I am grateful to Peter
Sarnak for calling my attention to this problem (a question
originally raised by Atle Selberg, cf. \cite{Se}), for his
comments about the paper and for many valuable discussions on
related topics. I also thank the referee for a careful reading and
for suggesting to clarify certain points in the argument.

\medskip \noindent {\bf Addendum (October 2001).} A preliminary
draft of this paper was completed in January 2001 and posted to
the \emph{e-Print archive} as 
\href{http://arxiv.org/abs/math.NT/0101096}{math.NT/0101096}. Later I
learned about the work of Kowalski, Michel and VanderKam
\cite{KMV} establishing subconvexity bounds for various families
of Rankin--Selberg $L$-functions. This work involves a more
elaborate version of the summation formula Proposition~\ref{Prop1}
below. The ultimate generalization (depending heavily on
Atkin--Lehner theory) appears in an unpublished complement
\cite{Mi} to \cite{KMV}. As pointed out on p.10 of \cite{Mi}, this
suffices, via the $\delta$-method, to establish Theorems~\ref{Th1}
and \ref{Th2} even in slightly stronger forms. In particular, the
original subconvexity exponent $1/2-1/22$ of \cite{DFI2} applies
for the general setting as well. However, the present paper is
technically simpler (e.g. it requires the theory of newforms only
to have the relevant $L$-functions at hand) and the simplification
is achieved by using Jutila's method of overlapping intervals in
place of the $\dd$-method.

\medskip \noindent {\bf Addendum (December 2002).} It is
straightforward to see from the argument given below that the
implied constants of Theorems 1 and 2 depend polynomially on $|s|$
and the levels of the forms involved. Some additional estimates on
Bessel functions establish polynomial dependence on the
Archimedean parameters (weight or Laplacian eigenvalue) as well.
The details are worked out for a special case in a recent paper by
Michel \cite{Mi2} where such a dependence turns out to be crucial.
\cite{Mi2} also supercedes the unpublished complement \cite{Mi}.

\section{Summation formula for the Fourier
coefficients}\label{sect1}

We define the normalized Fourier coefficients of cusp forms as
follows. Let $\phi(z)$ be a cusp form of level $N$ and nebentypus
$\chi$, that is, a holomorphic cusp form of some integral weight
$k$ or a real-analytic Maass cusp form of some nonnegative Laplacian
eigenvalue $1/4+\mu^2$. By definition, $\chi$ is a Dirichlet
character mod $N$ and the form $\phi$ satisfies a transformation
rule with respect to the Hecke congruence subgroup $\GG_0(N)$:
\[\phi\circ[\gG]=\chi(d)\phi,\qquad
\gG=\begin{pmatrix}a&b\\c&d\end{pmatrix}\in\GG_0(N),\]
where
\[\phi\circ[\gG](z)=\begin{cases}
\phi(\gG z)(cz+d)^{-k}&\text{if $\phi$ is holomorphic,}\\
\phi(\gG z)&\text{if $\phi$ is real-analytic,}
\end{cases}\]
and $\GG_0(N)$ acts on the upper half-plane $\HH=\{z:\Im z>0\}$ by
fractional linear transformations. Also, the form $\phi$ is
holomorphic or real-analytic on $\HH$ and decays exponentially to
zero at each cusp. Any such $\phi$ admits the Fourier expansion
\[\phi(z)=\sum_{m\neq 0}{\hat\lambda}_\phi(m)W(mz),\]
where
\[W(z)=\begin{cases}
e(z)&\text{if $\phi$ is holomorphic,}\\
|y|^{1/2}K_{i\mu}\bigl(2\pi|y|\bigr)e(x) &\text{if
$\phi$ is real-analytic.}
\end{cases}\]
Here $e(z)=e^{2\pi iz}$, $z=x+iy$ and $K_{i\mu}$ is the
MacDonald-Bessel function. If $\phi$ is holomorphic,
${\hat\lambda}_\phi(m)$ vanishes for $m\leq 0$. We define the
\emph{normalized Fourier coefficients} of $\phi$ as
\[\lfi(m)=\begin{cases}
{\hat\lambda}_\phi(m)m^\frac{1-k}{2}&\text{if $\phi$ is holomorphic,}\\
{\hat\lambda}_\phi(m)|m|^\frac{1}{2}&\text{if $\phi$ is
real-analytic.}
\end{cases}\]
This normalization corresponds to the Ramanujan Conjecture which
asserts that
\[\lfi(m)\ll_{\ep,\phi}m^\ep.\]
Rankin--Selberg theory implies that the conjecture holds on
average in the form
\begin{equation}\label{eq22}
\sum_{1\leq m\leq x}|\lfi(m)|^2\ll_\phi x.
\end{equation}

Various Voronoi-type summation formulas are fulfilled by these
coefficients. In the case of full level ($N=1$) Duke and Iwaniec
\cite{DI} established such a formula for holomorphic cusp forms
and Meurman \cite{M} for Maass cusp forms. These can be
generalized to arbitrary level and nebentypus with obvious minor
modifications as follows.
\begin{proposition}\label{Prop1}
Let $d$ and $q$ be coprime integers such that $N\mid q$, and let
$g$ be a smooth, compactly supported function on $(0,\infty)$. If
$\phi$ is a holomorphic cusp form of level $N$, nebentypus $\chi$
and integral weight $k$ then
\[\chi(d)\sum_{m=1}^\infty\lfi(m)e_q(dm)g(m)=
\sum_{m=1}^\infty\lfi(m)e_q\bigl(-{\bar d}m\bigr){\hat g}(m),\]
where
\[{\hat g}(y)=\frac{2\pi i^k}{q}\int_0^\infty g(x)
J_{k-1}\left(\frac{4\pi\sqrt{xy}}{q}\right)\,dx.\] If $\phi$ is a
real-analytic Maass cusp form of level $N$, nebentypus $\chi$ and
nonnegative Laplacian eigenvalue $1/4+\mu^2$ then
\[\chi(d)\sum_{m=1}^\infty\lfi(m)e_q(dm)g(m)=
\sum_{\pm}\sum_{m=1}^\infty\lfi(\mp m)e_q\bigl(\pm\bar d
m\bigr)g^{\pm}(m),\] where
\[\begin{split}
g^-(y)&=-\frac{\pi}{q\cosh \pi\mu}\int_0^\infty g(x)
\{Y_{2i\mu}+Y_{-2i\mu}\}\left(\frac{4\pi\sqrt{xy}}{q}\right)\,dx,\\\\
g^+(y)&=\frac{4\cosh \pi\mu}{q}\int_0^\infty g(x)
K_{2i\mu}\left(\frac{4\pi\sqrt{xy}}{q}\right)\,dx.
\end{split}\]
Here $\bar d$ is a multiplicative inverse of $d\bmod q$,
$e_q(x)=e(x/q)=e^{2\pi ix/q}$ and $J_{k-1}$, $Y_{\pm 2i\mu}$,
$K_{2i\mu}$ are Bessel functions.
\end{proposition}

The proof for the holomorphic case \cite{DI} is a straightforward
application of Laplace transforms. Meurman's proof for the
real-analytic case \cite{M} is more involved, but only because he
considers a wider class of test functions $g$ and has to deal with
delicate convergence issues. For smooth, compactly supported
functions $g$ as in our formulation these difficulties do not
arise and one can give a much simpler proof based on Mellin
transformation, the functional equations of the $L$-series
attached to additive twists of $\phi$ (see \cite{M}), and Barnes'
formulas for the gamma function. Indeed, Lemma 5 in \cite{St} (a
special case of Meurman's summation formula) has been proved by
such an approach. We expressed the formula for the non-holomorphic
case in terms of $K$ and $Y$ Bessel functions in order to
emphasize the analogy with the Voronoi-type formula for the
divisor function (where one has $\mu=0$) as derived by Jutila
\cite{Ju4,Ju5}.

\section{Setting up the circle method}\label{sect3}

For sake of exposition we shall only present the case of
real-analytic Maass forms and the equation $am-bn=h$. The other
cases follow along the same lines by changing Bessel functions and
signs at relevant places of the argument. In our inequalities
$\ep$ will always denote a small positive number whose actual
value is allowed to change at each occurrence. Furthermore, unless
otherwise indicated, implied constants will depend on $\ep$ and
the cusp forms only (including dependence on the level, nebentypus
characters and Laplacian eigenvalues).

Let $\phi(z)$ (resp. $\psi(z)$) be a Maass cusp form of level $N$,
nebentypus $\chi$ (resp. $\om$) and Laplacian eigenvalue
$1/4+\mu^2\geq 0$ (resp. $1/4+\nu^2\geq 0$) whose normalized
Fourier coefficients are $\lfi(m)$ (resp. $\lpsi(n)$), i.e.,
\[\begin{split}
\phi(x+iy)&=\sqrt{y}\sum_{m\neq 0}\lfi(m)K_{i\mu}\bigl(2\pi|m|y\bigr)e(mx),\\\\
\psi(x+iy)&=\sqrt{y}\sum_{n\neq
0}\lpsi(n)K_{i\nu}\bigl(2\pi|n|y\bigr)e(nx).
\end{split}\]

We shall first investigate $D_g(a,b;h)$ for smooth test functions
$g(x,y)$ which are supported in a box $[A,2A]\times[B,2B]$ and
have partial derivatives bounded by
\begin{equation}\label{eq27}g^{(i,j)}\ll_{i,j}A^{-i}B^{-j}P^{i+j}.
\end{equation} Our aim is to prove the estimate
\begin{equation}\label{eq29}D_g(a,b;h)\ll P^{11/10}(ab)^{-1/10}(A+B)^{1/10}(AB)^{2/5+\ep}.\end{equation}
In Section~\ref{sect2} we shall deduce Theorem~\ref{Th1} from this
bound by employing a partition of unity and decomposing
appropriately any smooth test function $f(x,y)$ satisfying
(\ref{eq1}). In fact, (\ref{eq29}) is a special case of
Theorem~\ref{Th1}, as can be seen upon setting $X=A$, $Y=B$, and
$f(x,y)=g(x,y)$. It supercedes the trivial upper bound
\[D_g(a,b;h)\ll(AB/ab)^{1/2}\] whenever
\begin{equation}\label{eq28}
ab\ll P^{-11/4}(A+B)^{-1/4}(AB)^{1/4-\ep}.
\end{equation}
The trivial bound is a consequence of $g\ll 1$, Cauchy's
inequality, and the Rankin--Selberg estimate (\ref{eq22}) applied
to the forms $\phi$ and $\psi$.

As $g(x,y)$ is supported in $[A,2A]\times[B,2B]$, we can assume
that $A,B\geq 1/2$, and also that
\begin{equation}\label{eq11}
h\leq 2(A+B),
\end{equation}
for otherwise $D_g(a,b;h)$ vanishes trivially. We shall attach, as
in \cite{DFI}, a redundant factor $w(x-y-h)$ to $g(x,y)$ where
$w(t)$ is a smooth function supported on $|t|<\dd^{-1}$ such that
$w(0)=1$ and $w^{(i)}\ll_i\dd^i$. This, of course, does not alter
$D_g(a,b;h)$. We choose
\begin{equation}\label{eq2}
\dd=P\frac{A+B}{AB},
\end{equation}
so that, by (\ref{eq27}), the new function
\[F(x,y)=g(x,y)w(x-y-h)\]
has partial derivatives bounded by
\begin{equation}\label{eq5}
F^{(i,j)}\ll_{i,j}\dd^{i+j}.
\end{equation}
We apply the Hardy--Littlewood method to detect the equation $am-bn=h$,
that is, we express $D_F(a,b;h)$ as the integral of a certain
exponential sum over the unit interval $[0,1]$. We get
\begin{equation}\label{eq30}D_g(a,b;h)=D_F(a,b;h)=\int_0^1
G(\al)\,d\al,\end{equation} where
\[G(\al)=\sum_{m,n}\lfi(m)\lpsi(n)F(am,bn)e\bigl((am-bn-h)\al\bigr).\]
We shall approximate this integral by the following proposition of Jutila
(a consequence of the main theorem in \cite{Ju3}).
\begin{proposition}[Jutila]\label{Prop2}
Let $\CQ$ be a nonempty set of integers $Q\leq q\leq 2Q$ where $Q\geq 1$.
Let $Q^{-2}\leq\dd\leq Q^{-1}$ and for
each fraction $d/q$ (in its lowest terms) denote by
$I_{d/q}(\al)$ the characteristic function of the interval
$\left[d/q-\dd,d/q+\dd\right]$. Write $L$ for the
number of such intervals, i.e.,
\[L=\sum_{q\in\CQ}\varphi(q),\]
and put
\[\tilde{I}(\al)=\frac{1}{2\dd L}\sum_{q\in\CQ}\ \
\csillag\sum_{d\smod{q}}I_{d/q}(\al).\] If $I(\al)$ is the
characteristic function of the unit interval $[0,1]$ then
\[\int_{-\infty}^\infty\bigl(I(\al)-\tilde{I}(\al)\bigr)^2 \,dx\ll
\dd^{-1}L^{-2}Q^{2+\ep},\]
where the implied constant depends on $\ep$ only.
\end{proposition}

We shall choose some $Q$ and apply the proposition with a set of
denominators of the form
\[\CQ=\bigl\{q\in[Q,2Q]:Nab\mid q\text{ and
}(h,q)=(h,Nab)\bigr\}.\] By a result of Jacobsthal \cite{Ja} the
largest gap between reduced residue classes mod $h$ is of size
$\ll h^\ep$, whence, by (\ref{eq11}),
\begin{equation}\label{eq18}
|\CQ|\gg\frac{Q(AB)^{-\ep}}{ab},
\end{equation}
assuming the right hand side exceeds some positive constant
$c=c(\ep,N)$. Moreover, we shall assume that
\begin{equation}\label{eq3}
Q^{-2}\leq\dd\leq Q^{-1},
\end{equation}
so that also
\begin{equation}\label{eq12}
1\leq Q\leq AB,
\end{equation}
whence (\ref{eq18}) yields
\begin{equation}\label{eq6}
L\gg\frac{Q^2(AB)^{-\ep}}{ab}.
\end{equation}
We clearly have
\begin{equation}\label{eq4}
|D_F(a,b;h)-\tilde{D}_F(a,b;h)|\leq\|G\|_\infty\|I-\tilde{I}\|_1,
\end{equation}
where
\[\begin{split}
\tilde{D}_F(a,b;h)&=\int_{-\infty}^\infty
G(\al)\tilde{I}(\al)\,d\al =\frac{1}{2\dd L}\sum_{q\in\CQ}\ \
\csillag\sum_{d\smod{q}}
\int_{-\infty}^\infty G(\al)I_{d/q}(\al)\,d\al\\\\
&=\frac{1}{2\dd L}\sum_{q\in\CQ}\ \ \csillag\sum_{d\smod{q}}
\int_{-\dd}^\dd G(d/q+\beta)\,d\beta =\frac{1}{2\dd
L}\sum_{q\in\CQ}\ \ \csillag\sum_{d\smod{q}}\FI_{d/q},
\end{split}\]
say. To derive an upper estimate for $G(\al)$ we express it as
\[G(\al)=\int_0^\infty\int_0^\infty F(x,y)e(-h\al)\,dS(x/a)\,dT(y/b),\]
where
\[S(x)=\sum_{1\leq m\leq x}\lfi(m)e(am\al),\quad
T(y)=\sum_{1\leq n\leq y}\lpsi(n)e(-bn\al).\] Then, integrating by
parts,
\[G(\al)=\int_0^\infty\int_0^\infty F^{(1,1)}(x,y)e(-h\al)S(x/a)T(y/b)\,dx\,dy,\]
therefore (\ref{eq5}) combined with Wilton's classical estimate
\[S(x)\ll x^{1/2}\log(2x),\quad T(y)\ll y^{1/2}\log(2y)\]
yields
\[\|G\|_\infty\ll\frac{(AB)^{1/2+\ep}}{(ab)^{1/2}}\|F^{(1,1)}\|_1
\ll\frac{\dd(AB)^{3/2+\ep}}{(ab)^{1/2}(A+B)}.\] Also, by
(\ref{eq6}) and Proposition~\ref{Prop2} we get
\[\|I-\tilde{I}\|_1\leq 3\|I-\tilde{I}\|_2\ll \frac{ab(AB)^\ep}{\dd^{1/2}Q},\]
so that (\ref{eq4}) becomes
\begin{equation}\label{eq15}
D_F(a,b;h)-\tilde{D}_F(a,b;h)\ll
\frac{(ab)^{1/2}\dd^{1/2}}{Q}\cdot\frac{(AB)^{3/2+\ep}}{A+B}.
\end{equation}

\section{Transforming exponential sums}

The contribution of the interval $[d/q-\dd,d/q+\dd]$ can be
expressed as
\[\FI_{d/q}=\int_{-\dd}^\dd
G(d/q+\beta)\,d\beta=e_q(-dh)\sum_{m,n}\lfi(m)\lpsi(n)e_q\bigl(d(am-bn)\bigr)E(m,n),\]
where
\[E(x,y)=F(ax,by)\int_{-\dd}^\dd
e\bigl((ax-by-h)\beta\bigr)\,d\beta.\] Using (\ref{eq5}) we
clearly have
\[E^{(i,j)}\ll_{i,j}\dd^{i+j+1}a^ib^j,\]
and we also record, for further reference, that
\begin{equation}\label{eq7}
\|E^{(i,j)}\|_1\ll_{i,j}\dd^{i+j}a^{i-1}b^{j-1}\frac{AB}{A+B}.
\end{equation}
We assume that $q\in\CQ$, hence $Nab\mid q$ and we can apply
Proposition~\ref{Prop1} to yield
\[\FI_{d/q}=\overline{\chi\om}(d)e_q(-dh)\sum_{\pm\pm}\,\sum_{m,n\geq
1}\lfi(\mp m)\lpsi(\mp n)e_q\bigl({\bar d}(\pm am\mp
bn)\bigr)E^{\pm\pm}(m,n),\] where the corresponding signs must be
matched and
\[E^{\pm\pm}(m,n)=\frac{ab}{q^2}\int_0^\infty\int_0^\infty
E(x,y)M^{\pm}_{2i\mu}\left(\frac{4\pi a\sqrt{mx}}{q}\right)
M^{\pm}_{2i\nu}\left(\frac{4\pi b\sqrt{ny}}{q}\right)\,dx\,dy,\]
\[M^+_{2ir}=(4\cosh\pi r)K_{2ir},\quad
M^-_{2ir}=-\frac{\pi}{\cosh\pi r}\{Y_{2ir}+Y_{-2ir}\}.\] By
summing over the residue classes we get
\begin{equation}\label{eq9}
\csillag\sum_{d\smod{q}}\FI_{d/q}= \sum_{\pm\pm}\,\sum_{m,n\geq
1}\lfi(\mp m)\lpsi(\mp n)S_{\overline{\chi\om}}(-h,\pm am\mp
bn;q)E^{\pm\pm}(m,n).
\end{equation}
In order to estimate the twisted Kloosterman sum we observe that
the greatest common divisor $(-h,\pm am\mp bn,q)$ divides
$N(h,n,a)(h,m,b)$ as follows from the relations $(a,b)=1$ and
$(h,q)=(h,Nab)$, therefore (\ref{eq10}) and (\ref{eq12}) imply
that
\begin{equation}\label{eq14}
S_{\overline{\chi\om}}(-h,\pm am\mp bn;q)\ll
(h,m)^{1/2}(h,n)^{1/2}Q^{1/2}(AB)^\ep.
\end{equation}
We estimate $E^{\pm\pm}(m,n)$ by successive applications of
integration by parts and the relations
\[\frac{d}{dz}\bigl(z^sK_s(z)\bigr)=-z^sK_{s-1}(z),
\quad\frac{d}{dz}\bigl(z^sY_s(z)\bigr)=z^sY_{s-1}(z);\]
\[K_s(z)\ll_s z^{-1/2},\quad Y_s(z)\ll_s z^{-1/2},\qquad z>0.\]
We get, for any integers $i,j\geq 0$,
\begin{multline*}E^{\pm\pm}(m,n)\ll_{i,j}
\frac{ab}{Q^2}\left(\frac{Q}{a\sqrt{m}}\right)^{i+\frac{1}{2}}
\left(\frac{Q}{b\sqrt{n}}\right)^{j+\frac{1}{2}}\\
\times\max_{\substack{0\leq k\leq i\\0\leq l\leq j}}
\left(\frac{A}{a}\right)^{k-\frac{i}{2}-\frac{1}{4}}
\left(\frac{B}{b}\right)^{l-\frac{j}{2}-\frac{1}{4}}
\|E^{(k,l)}\|_1,\end{multline*} i.e., by (\ref{eq7}),
\begin{multline*}E^{\pm\pm}(m,n)\ll_{i,j}\frac{1}{Q^2}
\left(\frac{Q}{a\sqrt{m}}\right)^{i+\frac{1}{2}}
\left(\frac{Q}{b\sqrt{n}}\right)^{j+\frac{1}{2}}\\
\times\left(\frac{A}{a}\right)^{-\frac{i}{2}-\frac{1}{4}}
\left(\frac{B}{b}\right)^{-\frac{j}{2}-\frac{1}{4}} \frac{AB}{A+B}
\max_{\substack{0\leq k\leq i\\0\leq l\leq
j}}(A\dd)^k(B\dd)^l.\end{multline*} Therefore
\begin{equation}\label{eq13}
E^{\pm\pm}(m,n)\ll_{i,j} \frac{(AB)^{1/2}}{\dd Q^2(A+B)}
\left(\frac{(\dd Q)^2A}{am}\right)^{\frac{i}{2}+\frac{1}{4}}
\left(\frac{(\dd Q)^2B}{bn}\right)^{\frac{j}{2}+\frac{1}{4}},
\end{equation}
suggesting that we can neglect the contribution to (\ref{eq9}) of
those pairs $(m,n)$ for which $am/A$ or $bn/B$ is $>(\dd
Q)^2(AB)^\ep$. Indeed, if we apply (\ref{eq22}) to $\phi$ and
$\psi$ to see that
\[\sum_{1\leq m\leq x}|\lfi(\mp m)|(h,m)^{1/2}\ll x\tau^{1/2}(h),\qquad
\sum_{1\leq n\leq y}|\lpsi(\mp n)|(h,n)^{1/2}\ll y\tau^{1/2}(h)\]
then we can specify $i$ and $j$ large enough (in terms of $\ep$)
to deduce from (\ref{eq14}) and (\ref{eq13}) that the contribution
to (\ref{eq9}) of those terms with $m$ or $n$ large is
\[\ll\tau(h)\frac{\dd^3Q^{5/2}(AB)^{-100}}{ab(A+B)},\]
say, while the choice $i=j=0$ in (\ref{eq13}) shows that the
remaining terms (for which $am/A$ and $bn/B$ are at most $(\dd
Q)^2(AB)^\ep$) contribute
\[\ll\tau(h)\frac{\dd^3Q^{5/2}(AB)^{3/2+\ep}}{ab(A+B)}.\]
Hence, by (\ref{eq6}),
\begin{equation}\label{eq16}
\tilde{D}_F(a,b;h)=\frac{1}{2\dd L}\sum_{q\in\CQ}\ \
\csillag\sum_{d\smod{q}}\FI_{d/q}
\ll\frac{\dd^2Q^{3/2}}{ab}\cdot\frac{(AB)^{3/2+\ep}}{A+B}.
\end{equation}

Inequalities (\ref{eq15}) and (\ref{eq16}) show that the optimal
balance is achieved when
\[\dd^3Q^5(ab)^{-3}\asymp 1.\]
A natural choice is given by
\[\dd^3Q^5=(cab)^3,\]
where $c$ is the constant appearing in the remark after
(\ref{eq18}). By (\ref{eq30}), this choice proves (\ref{eq29})
whenever the conditions of Proposition~\ref{Prop2} are satisfied,
that is, when $Q\geq cab(AB)^\ep$ and (\ref{eq3}) hold
simultaneously. It turns out that this is the case whenever
\[cab\leq P^{-2/3}(A+B)^{-2/3}(AB)^{2/3-\ep},\]
in particular, whenever (\ref{eq28}) is true. However, when
(\ref{eq28}) fails, (\ref{eq29}) follows from the Cauchy bound
$(AB/ab)^{1/2}$, as we already pointed out in Section~\ref{sect3}.

\section{Concluding Theorem~\ref{Th1}}\label{sect2}

Our aim is to prove Theorem~\ref{Th1} for all test functions
$f(x,y)$ satisfying (\ref{eq1}). We fix an arbitrary smooth
function
\[\rho:(0,\infty)\to\RR\]
whose support lies in $[1,2]$ and which satisfies the following
identity on the positive axis:
\[\sum_{k=-\infty}^\infty\rho\bigl(2^{-k/2}x\bigr)=1.\]
To obtain such a function, we take an arbitrary smooth
$\eta:(0,\infty)\to\RR$ which is constant 0 on $(0,1)$ and
constant 1 on $(\sqrt{2},\infty)$, and then define $\rho$ as
\[\rho(x)=\begin{cases}
\eta(x)&\text{if $0<x\leq\sqrt{2}$,}\\
1-\eta(x/\sqrt{2})&\text{if $\sqrt{2}<x<\infty$.}
\end{cases}\]
According to this partition of unity we decompose $f(x,y)$ as
\[f(x,y)=\sum_{k=-\infty}^\infty\sum_{l=-\infty}^\infty f_{k,l}(x,y),\]
\[f_{k,l}(x,y)=f(x,y)\rho\left(\frac{x}{2^{k/2}X}\right)
\rho\left(\frac{y}{2^{l/2}Y}\right).\] Observe that
\begin{equation}\label{eq8}\text{supp}\,f_{k,l}\subseteq
[A_k,2A_k]\times[B_l,2B_l\bigr],\quad A_k=2^{k/2}X,\quad
B_l=2^{l/2}Y,\end{equation} whence (\ref{eq1}) and $P\geq 1$ show
that
\[\bigl(1+2^{k/2}\bigr)\bigl(1+2^{l/2}\bigr)f_{k,l}^{(i,j)}\ll_{i,j}
A_k^{-i}B_l^{-j}P^{i+j}.\] In other words, the bound (\ref{eq29})
applies uniformly to each function
\[g_{k,l}(x,y)=\bigl(1+2^{k/2}\bigr)\bigl(1+2^{l/2}\bigr)f_{k,l}(x,y)\]
with the corresponding parameters $A=A_k$, $B=B_l$:
\[D_{g_{k,l}}(a,b;h)\ll
P^{11/10}(ab)^{-1/10}(A_k+B_l)^{1/10}(A_kB_l)^{2/5+\ep}.\] This
implies, for $\ep<1/10$,
\[D_{f_{k,l}}(a,b;h)\ll
2^{-|k|/5}2^{-|l|/5}P^{11/10}(ab)^{-1/10}(X+Y)^{1/10}(XY)^{2/5+\ep}.\]
Finally,
\[D_f(a,b;h)=\sum_{k=-\infty}^\infty\sum_{l=-\infty}^\infty
D_{f_{k,l}}(a,b;h)\] completes the proof of Theorem~\ref{Th1}.

It should be noted that the trivial upper bound
\[D_f(a,b;h)\ll(XY/ab)^{1/2}\] mentioned in Section~\ref{sect6}
follows by a similar reduction technique from the Cauchy bounds
$D_{g_{k,l}}(a,b;h)\ll(A_kB_l/ab)^{1/2}$ of Section~\ref{sect3}.

\section{Approximate functional equation}\label{sect4}

Let $\phi$ be a primitive holomorphic or Maass cusp
form of arbitrary level and nebentypus, $\Re s=1/2$, and $\chi$ a
primitive character modulo $q$ where $q$ is prime to the level.
Using the functional equation of the $L$-function attached to the
twisted primitive cusp form $\phi\otimes\chi$ and a standard
technique involving Mellin transforms we can express the special
value $L(s,\phi\otimes\chi)$ as a sum of two Dirichlet series of
essentially $\sqrt{C}$ terms where $C=C(s,\phi\otimes\chi)$ is the
analytic conductor defined by \cite{IS}. More precisely, $C\asymp
q^2$ where the implied constants depend only on $s$ and $\phi$,
therefore a special case of Theorem~2.1 in \cite{H} gives the
following
\begin{proposition}\label{Prop3}There is a smooth function
$f:(0,\infty)\to\CC$ and a complex number $\lL$ of modulus 1 such
that
\[L(s,\phi\otimes\chi)=\sum_{m=1}^\infty\frac{\lfi(m)\chi(m)}{m^{1/2}}
f\left(\frac{m}{q}\right)+
\lL\sum_{m=1}^\infty\frac{\blfi(m)\bar{\chi}(m)}{m^{1/2}}
\bar{f}\left(\frac{m}{q}\right).\] The function $f$ and its
partial derivatives $f^{(j)}$ $(j=1,2,.\dots)$ satisfy the
following uniform growth estimates at $0$ and infinity:
\[f(x)=\begin{cases}1+O(x^{\si}),&\quad
0<\si<1/5;\\
O(x^{-\si}),&\quad \si>0.
\end{cases}\]
\[f^{(j)}(x)\ll x^{-\si},\quad \si>j-1/5.\]
The implied constants depend only on $\si$, $j$, $s$ and the form
$\phi$.
\end{proposition}

For any positive numbers $A$ and $\ep$ we obtain, using
(\ref{eq22}), an expression
\[L(s,\phi\otimes\chi)=T+\lL\bar
T+O_{A,\ep,s,\phi}\left(q^{-A}\right),\] where
\[T=\sum_{1\leq m\leq q^{1+\ep}}
\frac{\lfi(m)\chi(m)g(m)}{m^{1/2}},\] and $g:(0,\infty)\to\CC$ is
a smooth function satisfying
\[g^{(j)}(x)\ll_{j,s,\phi} x^{-j}.\]
Therefore, applying partial summation and a smooth dyadic
decomposition we can reduce Theorem~\ref{Th2} to the following
\begin{proposition}\label{Prop4}
Let $1\leq M\leq q^{1+\ep}$ and $k$ be a smooth function supported
in $[M,2M]$ such that $k^{(j)}\ll_j M^{-j}$. Then
\[\sum_{m=1}^\infty\lfi(m)\chi(m)k(m)\ll q^{17/54+\ep}M^{2/3},\]
where the implied constant depends only on $\ep$ and the form
$\phi$.
\end{proposition}

\section{Amplification}

Our purpose is to prove Proposition~\ref{Prop4}. As in \cite{DFI2}
we shall estimate from both ways the amplified second moment
\[S=\csillag\sum_{\om\bmod{q}}
\left|\sum_{1\leq l\leq L}\bar\chi(l)\om(l)\right|^2|S_\om|^2,\]
where $\om$ runs through the primitive characters modulo $q$, $L$
is a parameter to be chosen later in terms of $M$ and $q$, and
\[S_\om=\sum_{m=1}^\infty\lfi(m)\om(m)k(m).\]
Assuming $L\geq c(\ep)q^\ep$ (indeed this will be the case
whenever $\ep<1/27$, cf. (\ref{eq26})) it follows, using the
result of Jacobsthal \cite{Ja} that the largest gap between
reduced residue classes mod $q$ is of size $\ll q^\ep$, that
\begin{equation}\label{eq20}
S\gg q^{-\ep}L^2|S_\chi|^2.
\end{equation}

On the other hand, expanding each primitive $\om$ in $S$ using
Gauss sums and then extending the resulting summation to all
characters mod $q$, we get by orthogonality,
\[S\leq\frac{\vfi(q)}{q}\csillag\sum_{d\smod{q}}
\left|\sum_n a(n)e_q(dn)\right|^2,\] where
\[a(n)=\sum_{\substack{lm=n\\1\leq l\leq
L}}\bar\chi(l)\lfi(m)k(m).\] It is clear that the coefficients
$a(n)$ are supported in the interval $[1,N]$ where $N=2LM$.
Extending the summation to all residue classes $d$ the previous
inequality becomes
\begin{equation}\label{eq21}
S\leq\vfi(q)\sum_{h\equiv 0\smod{q}}D(h),
\end{equation}
where
\[D(h)=\sum_{n_1-n_2=h}a(n_1)\bar a(n_2).\]

Using the Rankin--Selberg bound (\ref{eq22}) it is simple to
estimate the diagonal contribution $D(0)$. Indeed, by $k\ll 1$ we
get
\[\begin{split}
D(0)&=\sum_n|a(n)|^2\ll\sum_{\substack{l_1m_1=l_2m_2\\1\leq
l_1,l_2\leq L\\M\leq m_1,m_2\leq 2M}}\lfi(m_1)\blfi(m_2)\\\\
&\ll\sum_{\substack{1\leq l\leq L\\M\leq m\leq
2M}}|\lfi(m)|^2\tau(ml)\ll N^\ep L\sum_{M\leq m\leq
2M}|\lfi(m)|^2,
\end{split}\]
whence
\begin{equation}\label{eq23}
D(0)=\sum_n|a(n)|^2\ll N^{1+\ep}.
\end{equation}

In order to estimate the non-diagonal terms $D(h)$ $(h\neq 0)$ we
shall refer to Theorem~\ref{Th1}. Clearly, we can rewrite each
term as
\[D(h)=\sum_{1\leq l_1,l_2\leq L}\bar\chi(l_1)\chi(l_2)
\sum_{l_1m_1-l_2m_2=h}\lfi(m_1)\blfi(m_2)k(m_1)\bar k(m_2).\] The
inner sum is of type (\ref{eq24}), because $\blfi(m)$ is just the
$m$-th normalized Fourier coefficient of the contragradient cusp
form $\tilde\phi(z)=\bar\phi(-\bar z)$. For each pair $(l_1,l_2)$
we shall apply Theorem~\ref{Th1} with $a=l_1/(l_1,l_2)$,
$b=l_2/(l_1,l_2)$, $X=aM$ and $Y=bM$ to conclude that
\begin{equation}\label{eq25}
D(h)\ll L^2(a+b)^{1/10}(ab)^{3/10+\ep}M^{9/10+\ep}\ll
L^{27/10+\ep}M^{9/10+\ep}.
\end{equation}

\section{Concluding Theorem~\ref{Th2}}\label{sect5}
Inserting the bounds (\ref{eq23}) and (\ref{eq25}) into
(\ref{eq21}) it follows that
\[S\ll N^\ep\vfi(q)\left(N+\frac{N}{q}L^{27/10}M^{9/10}\right).\]
This shows that the optimal choice for $L$ is provided by
\begin{equation}\label{eq26}
q=L^{27/10}M^{9/10},
\end{equation}
whence (\ref{eq20}) yields
\[S_\chi\ll q^\ep L^{-1}|S|^{1/2}\ll
(qN)^{1/2+\ep}L^{-1}\ll (qM/L)^{1/2+\ep}.\] Substituting
(\ref{eq26}) we get
\[S_\chi\ll (qMq^{-10/27}M^{1/3})^{1/2+\ep}\ll q^{17/54+\ep}M^{2/3},\]
which is precisely the conclusion of Proposition~\ref{Prop4}. The
proof of Theorem~\ref{Th2} is complete.

\end{document}